\newcommand{\papertitle}{A log-concavity criterion for finite nonnegative linearization}
\newcommand{\papershort}{A LOG-CONCAVITY CRITERION FOR FINITE LINEARIZATION}
\newcommand{\paperdate}{November 1, 2026}
\documentclass[twoside,10pt]{article}
\usepackage[T1]{fontenc}
\usepackage{lmodern}
\usepackage{amsmath,amssymb,amsthm,amsfonts,mathtools}
\usepackage[letterpaper,total={5.5in,8in},vmarginratio=1:1,hmarginratio=1:1,headsep=21pt]{geometry}
\usepackage{xcolor}
\usepackage{array,longtable}
\usepackage{xr-hyper}
\usepackage[breaklinks=true,linktocpage=true]{hyperref}
\hypersetup{colorlinks=true,linkcolor=blue!45!black,urlcolor=blue!45!black,citecolor=blue!45!black,
pdfauthor={Leslie P. Polzer},pdftitle={\papertitle},
pdfsubject={Finite orthogonal polynomials and nonnegative linearization},
pdfkeywords={orthogonal polynomials, Krawtchouk polynomials, nonnegative linearization, quadratic transformation, hypergroups, Lean}}
\numberwithin{equation}{section}
\newtheoremstyle{slbody}{8pt}{8pt}{\normalfont\slshape}{}{\bfseries}{.}{0.5em}
{\thmname{#1}\thmnumber{ #2}\thmnote{ {\normalfont\itshape(\ignorespaces#3\unskip)}}}
\theoremstyle{slbody}
\newtheorem{theorem}{Theorem}[section]
\newtheorem{lemma}[theorem]{Lemma}
\newtheorem{corollary}[theorem]{Corollary}
\newtheorem{proposition}[theorem]{Proposition}
\theoremstyle{definition}

\theoremstyle{remark}
\newtheorem{remark}[theorem]{Remark}

\newcommand{\xnum}[2]{#2}
\title{\papertitle}
\author{Leslie P. Polzer}
\newcommand{\authoraffiliation}{Independent Researcher}
\newcommand{\authoremail}{polzer@fastmail.com}
\date{\paperdate}
\makeatletter
\renewcommand{\section}{\@startsection{section}{1}{\z@}{-3.5ex \@plus -1ex \@minus -.2ex}{2.3ex \@plus .2ex}{\normalfont\large\bfseries}}
\renewcommand{\@seccntformat}[1]{\csname the#1\endcsname.\enspace}
\def\tagform@#1{\maketag@@@{(\ignorespaces{\oldstylenums{#1}}\unskip\@@italiccorr)}}
\renewcommand{\eqref}[1]{\textup{{\normalfont(\oldstylenums{\ref{#1}}}\normalfont)}}
\newcommand{\shorttitle}{\papershort}
\newcommand{\shortauthors}{LESLIE P. POLZER}
\newcommand{\ps@bookheader}{%
\renewcommand\@oddfoot{\hfil}\renewcommand\@evenfoot{\hfil}%
\renewcommand\@oddhead{\ifnum\value{page}>1{\small\hfil\shortauthors}\hfil\thepage\else\hfil\fi}%
\renewcommand\@evenhead{\thepage\hfil\small\shorttitle\hfil}}
\renewcommand{\maketitle}{\begin{center}{\large\bfseries\@title}\\\vskip36pt
{\scshape\@author}\par\medskip{\small\authoraffiliation}\par\smallskip
{\small\texttt{\authoremail}}\par\smallskip{\small\@date}\par\vskip36pt\end{center}}
\makeatother
\renewenvironment{abstract}{\quotation\noindent\small{\bfseries\abstractname.\enspace}}{\endquotation}
\usepackage{tikz}
\theoremstyle{slbody}
\newtheorem*{mainthm}{Theorem~A}
\begin{document}
\maketitle

\begin{abstract}
In a 2005 survey, R. Szwarc asked for criteria, stated in the recurrence coefficients, for nonnegative linearization (property (P)) of finite orthogonal polynomial systems such as the Krawtchouk polynomials. For a finite three-term recurrence normalized to have row sums one, whose diagonal and downward coefficients have a constant nonnegative ratio, we prove that property (P) holds for every product, including products beyond the terminal degree, whenever the Newton coefficients of the profile $d_i=b_i+c_i$ are log-concave without internal zeros. The criterion contains the full positive Krawtchouk range and the dual $q$-Krawtchouk arrays, and through a quadratic decomposition of symmetric finite systems it covers Szwarc's quadratic system $Q_n$ in its sharp range $0<p_0\leq1/2$. As a limit to any such criterion, property (P) fails at coefficient distance $O((1-r^2)^N)$ from the Krawtchouk array of size $N$ with fixed ratio $c_i/d_i=r<1$. All theorems, corollaries and propositions are formalized in Lean~4.
\end{abstract}

\vskip50pt
\baselineskip=13.5pt
\section{Introduction}\label{sec:introduction}

Nonnegative linearization (property (P)) asks whether the product of two orthogonal polynomials is a nonnegative linear combination of polynomials of the same family. For a finite system $R_0,\ldots,R_N$ the identity is understood on the orthogonality support, so it includes products whose degree sum exceeds $N$. In the normalization $R_i(1)=1$ the coefficients then sum to one and define a finite polynomial hypergroup.

In his survey, Szwarc~\cite[\S5, Problem~2]{szwarc} asks: ``Find any criteria for the property (P) for finite systems of orthogonal polynomials, like the Krawtchouk polynomials.'' He then introduces a symmetric system $Q_n$ built from the Krawtchouk recurrence, notes that it has property (P) for $p_0\leq1/2$, and asks for a general criterion that would handle this case.

The main result answers this for arrays whose diagonal and downward coefficients have a constant nonnegative ratio. For a normalized recurrence \eqref{eq:recurrence}, whose coefficients have row sums one, write $d_i=b_i+c_i$ for the profile and $\nu_k=\Delta^kd(0)=\sum_{i\leq k}(-1)^{k-i}\binom kid_i$ for its Newton coefficients.

\begin{mainthm}[Theorem~\ref{thm:lc}]
Let $0=d_0<d_1<\cdots<d_N=1$ have Newton coefficients that are nonnegative and log-concave without internal zeros, let $0<r\leq1$, and put $a_i=1-d_i$, $b_i=(1-r)d_i$ and $c_i=rd_i$. Then every product $R_mR_n$, including those with $m+n>N$, is a nonnegative combination of $R_0,\ldots,R_N$ on the orthogonality support.
\end{mainthm}

The input is the profile alone: one computes the $N$ numbers $\nu_k$ and checks their signs and $N-2$ quadratic inequalities, without forming any product coefficient. For instance, $d=\frac1{97}(0,8,24,52,97)$ has $\nu=\frac1{97}(8,8,4,1)$, which is log-concave, so all its constant-ratio arrays have property (P); this profile is neither Krawtchouk, nor geometric, nor quadratic, and Section~\ref{sec:lc} embeds it in a one-parameter class of such profiles.

Property (P) of the classical examples is known; what is new is a coefficient criterion that recovers it, and a reduction for symmetric finite systems. Theorem~A contains the full positive Krawtchouk range (Corollary~\ref{cor:kraw}), whose property (P) is classical~\cite{eagleson,dunkl,szwarc}, and the geometric arrays of Corollary~\ref{cor:additional}, the dual $q$-Krawtchouk systems of Koelink's construction~\cite{koelink}. For the second part of Szwarc's problem, a quadratic decomposition reduces property (P) of a symmetric finite system to that of two nonsymmetric ones (Theorem~\ref{thm:decomp}); with Theorem~A it gives a parity-split criterion (Corollary~\ref{cor:parity}) that covers $Q_n$ in its sharp range $0<p_0\leq1/2$ (Theorem~\ref{thm:szwarc}), whose positivity Szwarc states. Table~\ref{tab:compare} compares Theorem~A with earlier criteria.

\begin{table}[!t]
\centering
\caption{Theorem~A and earlier criteria for nonnegative linearization, as they apply to finite systems.}\label{tab:compare}
\renewcommand{\arraystretch}{1.15}
\begin{tabular}{@{}>{\raggedright\arraybackslash}p{0.22\linewidth}>{\raggedright\arraybackslash}p{0.38\linewidth}>{\raggedright\arraybackslash}p{0.32\linewidth}@{}}
Criterion & Hypothesis & Products of a finite system\\\hline
Theorem~A & inequalities on the Newton coefficients of the profile; constant ratio $b_i/c_i$ & all, including degree sum $>N$\\
Askey~\cite[Part~II]{askey} & conditions on the recurrence coefficients, through a maximum principle & degree sum $\leq N$\\
Szwarc~\cite[Thms.~3.2, 3.3]{szwarc} & monotonicity conditions for infinite systems & not applicable: some conditions fail at the terminal index\\
Szwarc~\cite[Thm.~1]{szwarc2003} & an auxiliary discrete hyperbolic problem; necessary and sufficient & all; not an inequality test on the coefficients\\
M{\l}otkowski--Szwarc~\cite{mlotszwarc} & chain-sequence conditions for infinite systems on convergent point sequences & not applicable (infinite systems)\\
Miclo~\cite[Prop.~3.3]{miclo} & spectral; necessary and sufficient & all; not an inequality test on the coefficients\\
Miclo~\cite[Thm.~31]{miclo2019} & inequality conditions for chains symmetric under $i\mapsto N-i$ & all; excludes the Krawtchouk systems with $N>2$
\end{tabular}
\end{table}

\emph{Related work.} Szwarc's monotonicity criteria fail for finite systems near the terminal index, whose upward coefficient vanishes (Section~\ref{sec:szwarc}). Miclo, having found numerically that two natural discrete versions of the log-concavity theorem of Achour and Trim\`eche fail, concluded that ``the right notion of log-concavity of a measure (or rather of a Markov kernel) has yet to be found in the discrete setting''~\cite[p.~433]{miclo}; his condition in~\cite[Theorem~31]{miclo2019} that $2^ia_i$ be nonincreasing on $[0,N/2]$ fails for the symmetric Krawtchouk system once $N>2$, and the systems with $p_0\neq1/2$ are not symmetric. For the finite systems of the Askey scheme a related question is open: Terwilliger~\cite[Problem~36.4]{terwilliger} asks which real Leonard systems have only nonnegative intersection numbers. Blau and Hein~\cite[Theorem~1.1]{blauhein} determine the recurrences of the homogeneous monotonic $P$-polynomial table algebras, which are finite systems with nonnegative linearization in which all nontrivial elements have the same Haar weight and the recurrence coefficients are monotone.

Section~\ref{sec:limits} marks the scope of the criterion: for fixed $r<1$, coefficient balls around $K_{N,r}$ in which every array has property (P) have radius at most $3(1-r^2)^N$ for large $N$ (Appendix~\ref{app:radius}), so a criterion valid for all $N$ must place these arrays near its boundary, and Theorem~A places them on its edge. For fixed $N$ and $r\to1$ the order is $(1-r)^N$~\cite[Corollary~\xnum{P3-cor:radius-order}{9.11}]{P3}. Positivity of the Newton coefficients does not suffice, but log-concavity is not necessary; we pose the exact boundary as a problem.

\emph{Companion papers and formalization.} The written proofs here do not use the companion papers, which are cited for comparison only. The companion paper~\cite{P2} builds on this one and extends Theorem~A to ratios $c_i/d_i$ that increase with the index; the companion paper~\cite{P3} is complementary: it gives explicit perturbation neighborhoods, inequality and reversal criteria, and constraints at the symmetric Krawtchouk endpoint. All theorems, corollaries and propositions are formalized in Lean~4; the section on code and formalization states which Lean proofs take routes different from the written ones and how the three papers depend on each other.

\section{Setting}\label{sec:criterion}

For integers $j\geq0$ put $[j]_q=1+q+\cdots+q^{j-1}$, with $[0]_q=0$; in particular $[j]_1=j$.
Consider real coefficients indexed by $0\leq i\leq N$, where $N\geq1$, satisfying
\begin{equation}\label{eq:basic}
\begin{gathered}
a_0=1,\quad b_0=c_0=0,\quad a_N=0,\quad a_i+b_i+c_i=1,\\
a_i>0\ (i<N),\qquad c_i>0\ (i>0).
\end{gathered}
\end{equation}
Set $R_{-1}=0$, $R_0=1$, and determine $R_1,\ldots,R_N$ by
\begin{equation}\label{eq:recurrence}
xR_i=a_iR_{i+1}+b_iR_i+c_iR_{i-1}\qquad(0\leq i<N).
\end{equation}
At $i=0$, \eqref{eq:recurrence} gives $R_1=x$; induction at $x=1$, using the row sums and $a_i>0$ for $i<N$, gives $R_i(1)=1$ for $0\leq i\leq N$.
The terminal relation is $(x-b_N)R_N-c_NR_{N-1}=0$ on the orthogonality support. Linearization means equality of functions on this support, rather than an unrestricted polynomial identity when the degree sum exceeds $N$.

Write the tridiagonal recurrence matrix of such an array with the row convention
\[
J_{i,i+1}=a_i,\qquad J_{i,i}=b_i,\qquad J_{i,i-1}=c_i,
\]
with all other entries zero and indices outside $0,\ldots,N$ omitted. We call an array satisfying \eqref{eq:basic} \emph{normalized}. Its matrix $J$ has row sums one, but \eqref{eq:basic} allows $b_i<0$; we call the array \emph{stochastic} only when all $b_i\geq0$, so that $J$ is a stochastic matrix. The arrays of Theorem~\ref{thm:lc} below are stochastic.

\begin{lemma}[Finite recurrence realization and triangle support]\label{lem:finite}
Every recurrence satisfying \eqref{eq:basic} has a positive orthogonality measure on exactly $N+1$ real points. If $J$ is its tridiagonal recurrence matrix, then its finite product coefficients satisfy
\begin{equation}\label{eq:matrix-product}
g(m,n;k)=[R_m(J)]_{n,k}.
\end{equation}
They vanish unless $|m-n|\leq k\leq\min(m+n,N)$, regardless of their signs inside this interval.
\end{lemma}
\begin{proof}
Choose positive $\pi_i$ with $\pi_0=1$ and $\pi_{i+1}=\pi_i a_i/c_{i+1}$. For $D=\operatorname{diag}(\sqrt{\pi_i})$, $H=DJD^{-1}$ is real symmetric. Induction in \eqref{eq:recurrence} gives the row-vector identity $e_0^\mathsf T R_i(J)=e_i^\mathsf T$, hence $e_0^\mathsf T R_i(H)=\pi_i^{-1/2}e_i^\mathsf T$. In particular $e_0$ is cyclic for $H$. In an orthogonal eigenspace decomposition of $H$, the cyclic span of $e_0$ has dimension at most the number of distinct eigenvalues; here it is all of $\mathbb R^{N+1}$, so $H$ has $N+1$ distinct eigenvalues and the projection of $e_0$ on each eigenspace is nonzero. The spectral measure $\mu$ of $H$ at $e_0$ therefore has positive mass at each of $N+1$ real points, and
\[
\int R_iR_j\,d\mu=e_0^\mathsf TR_i(H)R_j(H)e_0
=\boldsymbol1_{i=j}\pi_i^{-1}.
\]
The terminal relation holds on this support: its polynomial $Q(x)=(x-b_N)R_N-c_NR_{N-1}$ satisfies $e_0^\mathsf TQ(J)=0$. Since $e_0$ is cyclic and $Q(J)$ commutes with $J$, $Q(J)=0$.
Multiplying a basis polynomial by $x$ acts by the corresponding row of $J$, including the terminal row. Induction on powers of $x$ therefore gives \eqref{eq:matrix-product}. Orthogonality also gives
\[
g(m,n;k)=\frac{\int R_mR_nR_k\,d\mu}{\int R_k^2\,d\mu}.
\]
If one of $m,n,k$ exceeds the sum of the other two, its polynomial is orthogonal to their product, whose degree is smaller. This proves the asserted support.
\end{proof}

We say that the array has \emph{nonnegative linearization}, or \emph{property (P)}, if $g(m,n;k)\geq0$ for all $0\leq m,n,k\leq N$. The coefficients $g(m,n;\cdot)$ always sum to one, since $R_i(1)=1$ and $1$ lies in the support, because $J\boldsymbol1=\boldsymbol1$.

For $q\geq1$ and $0<r\leq1$ let $G_{N,q,r}$ be the array
\begin{equation}\label{eq:parameters}
a_i=\frac{q^i[N-i]_q}{[N]_q},\qquad
b_i=(1-r)\frac{[i]_q}{[N]_q},\qquad
c_i=r\frac{[i]_q}{[N]_q},
\end{equation}
with $[j]_q$ as above; it satisfies \eqref{eq:basic}, since $[N]_q-[i]_q=q^i[N-i]_q$. Its profile $d_i=[i]_q/[N]_q$ has geometric increments with ratio $q$, and $b_i/c_i=(1-r)/r$ is constant. The Krawtchouk array is $K_{N,r}=G_{N,1,r}$, with $a_i=(N-i)/N$, $b_i=(1-r)i/N$ and $c_i=ri/N$.

\section{Log-concave profiles}\label{sec:lc}

For an array satisfying \eqref{eq:basic} with profile $d_i=b_i+c_i$ write
\[
\nu_k=\sum_{i=0}^k(-1)^{k-i}\binom ki d_i\quad(1\leq k\leq N),\qquad d_n=\sum_{k=1}^n\binom nk\nu_k\quad(0\leq n\leq N);
\]
below we put $\nu_0=0$ and $\nu_k=0$ for $k>N$.
We say that $(\nu_k)_{1\leq k\leq N}$ is \emph{log-concave without internal zeros} if all $\nu_k\geq0$, the indices with $\nu_k>0$ form an interval $\{1,\ldots,K\}$, and $\nu_k^2\geq\nu_{k-1}\nu_{k+1}$ for $2\leq k\leq N-1$. Note $\nu_1=d_1>0$.

\begin{theorem}[Log-concave Newton criterion]\label{thm:lc}
Let $0=d_0<d_1<\cdots<d_N=1$ have Newton coefficients that are log-concave without internal zeros, let $0<r\leq1$, and put $a_i=1-d_i$, $b_i=(1-r)d_i$, $c_i=rd_i$. Then the system has nonnegative linearization.
\end{theorem}

The hypotheses are inequalities in the profile. They hold for the Krawtchouk arrays ($\nu=(1/N,0,\ldots,0)$) and for all arrays $G_{N,q,r}$ ($\nu_k$ proportional to $(q-1)^{k-1}$, log-concave with equality). The hypotheses also hold for every admissible convex quadratic profile: $\nu_1>0$, $\nu_2\geq0$ when $N\geq2$, and $\nu_k=0$ for $3\leq k\leq N$. The hypotheses exclude the known counterexamples with constant ratio: the convex counterexamples recorded in~\cite{P3} have Newton coefficients of both signs, and for the combinations $(1-\lambda)K_{N,r}+\lambda G_{N,q,r}$ of Section~\ref{sec:limits}, with $q>1$, $N\geq3$ and $0<\lambda<1$, one has $\nu_1=(1-\lambda)/N+\lambda/[N]_q$ and $\nu_k=\lambda(q-1)^{k-1}/[N]_q$ for $2\leq k\leq N$, so $\nu_2^2-\nu_1\nu_3=-\lambda(1-\lambda)(q-1)^2/(N[N]_q)<0$.

\medskip
\noindent\fbox{\begin{minipage}{0.965\linewidth}
\textbf{Applying Theorem~\ref{thm:lc}.} Given a profile $0=d_0<d_1<\cdots<d_N=1$:
\begin{enumerate}
\item compute $\nu_k=\sum_{i\leq k}(-1)^{k-i}\binom kid_i$ for $1\leq k\leq N$;
\item check that all $\nu_k\geq0$, that the $k$ with $\nu_k>0$ form an interval $\{1,\ldots,K\}$, and that $\nu_k^2\geq\nu_{k-1}\nu_{k+1}$ for $2\leq k\leq N-1$;
\item conclude property (P), for all products $R_mR_n$ with $0\leq m,n\leq N$, for $a_i=1-d_i$, $b_i=(1-r)d_i$, $c_i=rd_i$ and every $0<r\leq1$.
\end{enumerate}
For $d=\frac1{97}(0,8,24,52,97)$: $97\nu_1=8$, $97\nu_2=24-16=8$, $97\nu_3=52-72+24=4$, $97\nu_4=97-208+144-32=1$. All are positive, and $8^2\geq8\cdot4$, $4^2\geq8\cdot1$. So every constant-ratio array with this profile has property (P); at $r=1$, for instance, $g(4,4;4)=95753/292365$.
\end{minipage}}
\medskip

This profile is a member of a one-parameter class. For $\alpha>0$ and $0<q\leq1$ let $\nu_k=C\alpha^{k-1}q^{\binom{k-1}2}$ for $1\leq k\leq N$, with $C>0$ chosen so that $d_N=\sum_k\binom Nk\nu_k=1$. Then $\nu_k^2=q^{-1}\nu_{k-1}\nu_{k+1}\geq\nu_{k-1}\nu_{k+1}$, so Theorem~\ref{thm:lc} applies. For $q=1$ these are the profiles of the geometric arrays $G_{N,1+\alpha,r}$, and as $q\to0$ they tend to quadratic profiles. For $0<q<1$ and $N\geq3$ the increments are not geometric and $\nu_3\neq0$, so the arrays are neither Krawtchouk, nor geometric, nor quadratic. The profile in the box is $N=4$, $\alpha=1$, $q=\frac12$.

\begin{figure}[t]
\centering
\begin{tikzpicture}[x=0.95cm,y=0.62cm,>=stealth]
\foreach \l in {0,...,7} \foreach \a in {0,...,7} {\pgfmathtruncatemacro{\s}{\l+\a}\ifnum\s<8 \fill[black!30] (\l,\a) circle (1.3pt);\fi}
\draw[black!45] (0,7.5) -- (0,0) -- (7.5,0);
\draw[black!30] (0,7) -- (7,0);
\node[anchor=north] at (7.5,0) {$\ell$};
\node[anchor=east] at (0,7.5) {$a$};
\node[black!55,anchor=south west,font=\small] at (0.05,7.05) {$\ell+a=N$};
\fill (1,6) circle (2.2pt) node[anchor=east,font=\small] {$x$\ };
\fill (2,4) circle (2.2pt) node[anchor=east,font=\small] {$z$\ };
\fill (3,2) circle (2.2pt) node[anchor=north west,font=\small] {$y$};
\fill (2,2) circle (2.2pt);
\draw[->,thick] (2,2) -- (2.9,2);
\node[anchor=north,font=\small] at (2.45,1.95) {$\mathcal U$};
\draw[->,thick] (1,6) -- (1.95,4.1);
\draw[->,thick] (2,4) -- (2.95,2.1);
\node[anchor=east,font=\small] at (1.4,5.0) {$B_0$};
\node[anchor=east,font=\small] at (2.4,3.0) {$B_0$};
\draw[->,thick,dashed] (1,6) to[bend left=22] (2.93,2.12);
\node[anchor=west,font=\small] at (2.45,4.4) {$-B_1$};
\end{tikzpicture}
\caption{The kernel recurrence of Lemma~\ref{lem:kernel-rec} on the triangle $\Omega$ ($N=7$). Every move raises $\ell$. The node $y=(\ell+1,a)$ receives the straight move $\mathcal U$ from $(\ell,a)$ and moves of index $j$ from nodes with smaller $\ell$; the moves of odd index carry a minus sign. The dashed move $-B_1$ from $x$, of displacement $(2,-4)$, is dominated by a quarter of the two-step paths from $x$ through $z=x+(1,-2)$ to $y$, made of two moves $B_0$ of even index (Lemma~\ref{lem:detours}); this gives $\mathcal P^2\geq4\mathcal Q$.}\label{fig:moves}
\end{figure}
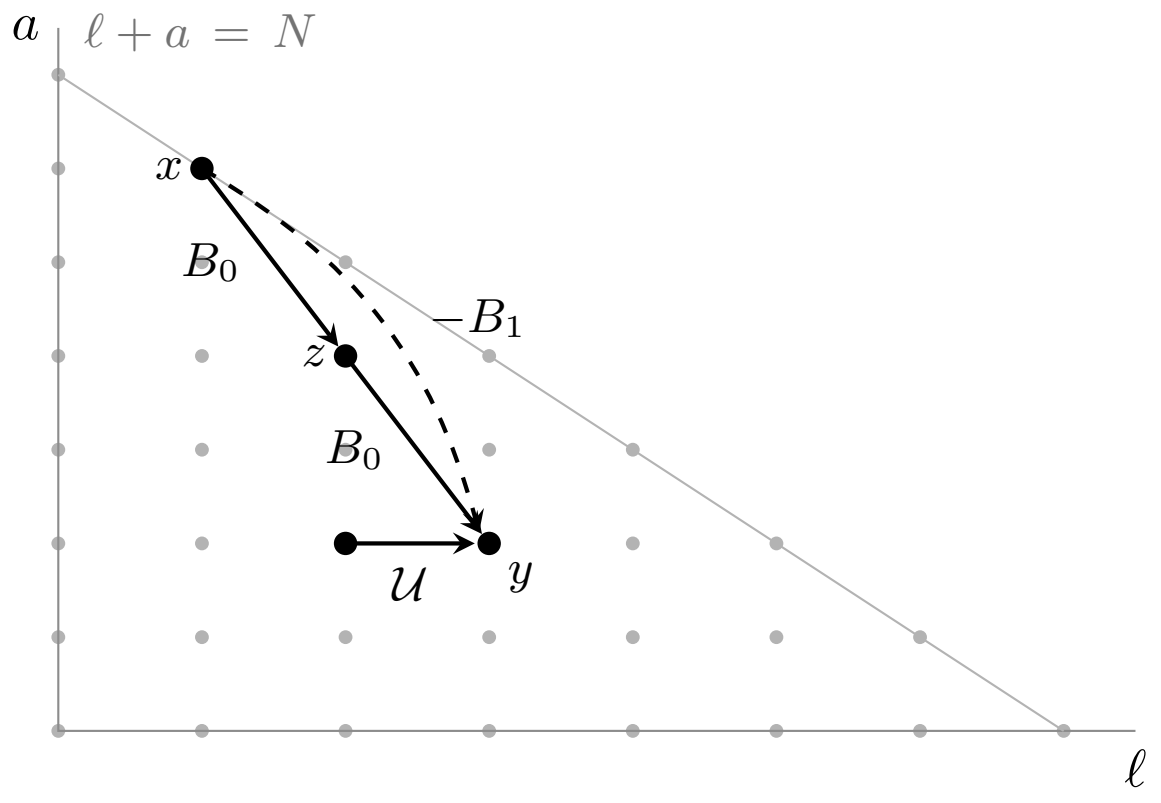

\paragraph{Outline of the proof.}
The proof has four steps. (1) A rescaling removes $r$ from the downward coefficients, and in the basis $E_i=e_ip_i$ the product coefficients become positive multiples of numbers $L_{mn}^k$. For fixed $k$ these form a symmetric generating function $K_k(z,w)$, which an explicit difference operator annihilates; the operator includes the terminal row, so the argument also covers products of degree sum greater than $N$. Expanding $K_k=\sum\omega_{\ell,a}(z+w)^a(zw)^\ell$, it suffices to show $\omega_{\ell,a}\geq0$. (2) Lemma~\ref{lem:kernel-rec} turns the annihilation into a recurrence that expresses $\omega$ at level $\ell+1$ through lower levels, by \emph{moves} with explicit weights $a_{\ell+a}$, $tA_j$, $B_j$ and $C_j$; the moves of odd index $j$ carry a minus sign. The recurrence holds for every profile, and its weights are nonnegative when all $\nu_k\geq0$. (3) Log-concavity enters only through Lemmas~\ref{lem:lc-ineq} and~\ref{lem:detours}: each move of odd index $j$ is dominated by a quarter of the two-step paths through moves of the even indices $0$ and $j-1$. (4) In operator form, $\omega=\beta+(\mathcal U+\mathcal P-\mathcal Q)\omega$ with $\mathcal P^2\geq4\mathcal Q$, and induction on $\ell$ gives $\omega\geq\frac12\mathcal P\omega\geq0$. Figure~\ref{fig:moves} shows the moves and one detour, and Table~\ref{tab:detours} lists all detours.

\begin{table}[t]
\centering
\caption{The negative moves into $y=(\ell+1,a)$, for odd $j$, and the two-step paths through moves of even index that dominate them (Lemma~\ref{lem:detours}). The first step of a detour has index $0$, the second has index $j-1$.}\label{tab:detours}
\renewcommand{\arraystretch}{1.15}
\begin{tabular}{@{}lll@{}}
Negative move & Source node & Dominating detours\\\hline
$tA_j$ & $(\ell-j,\,a+2j+1)$ & $A_0$ then $B_{j-1}$; \ $B_0$ then $A_{j-1}$\\
$B_j$ & $(\ell-j,\,a+2j+2)$ & $B_0$ then $B_{j-1}$\\
$C_j$ & $(\ell-j-1,\,a+2j+2)$ & $B_0$ then $C_{j-1}$
\end{tabular}
\end{table}

\paragraph{A rescaled system and a kernel.}
Fix $r$ and put $t=(1-r)/\sqrt r\geq0$. If $p_i$ is defined by $p_0=1$ and $yp_i=a_ip_{i+1}+td_ip_i+d_ip_{i-1}$ $(0\leq i<N)$, then $R_i(x)=r^{i/2}p_i(x/\sqrt r)$, so the coefficients $H_{mn}^k$ in $p_mp_n=\sum_kH_{mn}^kp_k$ on the support satisfy $g(m,n;k)=r^{(m+n-k)/2}H_{mn}^k$. It suffices to show $H_{mn}^k\geq0$. Put
\[
e_i=\prod_{j=0}^{i-1}\frac{a_j}{d_{j+1}},\qquad E_i=e_ip_i.
\]
Then $yE_i=d_{i+1}E_{i+1}+td_iE_i+a_{i-1}E_{i-1}$ on the support for $0\leq i\leq N$, with the terms $d_{N+1}E_{N+1}$ and $a_{-1}E_{-1}$ omitted; at $i=N$ this is the terminal relation. Write $E_mE_n=\sum_kL_{mn}^kE_k$, so that $H_{mn}^k=L_{mn}^ke_k/(e_me_n)$, and
\[
K_k(z,w)=\sum_{m,n=0}^NL_{mn}^kz^mw^n.
\]
Let $\theta_z=z\,\partial_z$ and let $d(\theta)=\sum_k\nu_k\binom\theta k$ act on polynomials by $d(\theta_z)z^i=d_iz^i$. For $F(z)=\sum_{i=0}^NE_iz^i$ the recurrence reads $yF(z)=\mathcal L_zF(z)$ with
\[
\mathcal L_z=(z^{-1}+t-z)\,d(\theta_z)+z,
\]
because $d_0=0$ removes the term $z^{-1}E_0$ and $1-d_N=0$ removes the term $z^{N+1}E_N$. Applying this to $F(z)F(w)=\sum_kK_k(z,w)E_k$ in the two variables separately and comparing coefficients of the basis $E_k$ gives $(\mathcal L_z-\mathcal L_w)K_k=0$. Since $K_k$ is symmetric, $K_k=\sum_{\ell,a\geq0}\omega_{\ell,a}S^aZ^\ell$ with $S=z+w$, $Z=zw$, and $K_k(z,0)=z^k$ gives $\omega_{0,a}=\mathbf1_{a=k}$. As $K_k$ has degree at most $N$ in $z$, comparing the coefficients of $z^{\ell+a}w^\ell$ by induction on $\ell$ gives $\omega_{\ell,a}=0$ whenever $\ell+a>N$.

\paragraph{The kernel recurrence.}
Put $f_n(\ell)=\sum_{h\geq0}\binom\ell h\nu_{n+h}$ with $\nu_0=0$ and $\nu_k=0$ for $k>N$, so $f_0(\ell)=d_\ell$, and
\[
W_{p,q}(\ell,a)=\sum_{n\geq p}\binom{n-p+q}q\binom an f_n(\ell)\geq0,
\]
and $A_j=W_{2j+1,j}$, $B_j=W_{2j+2,j}$, $C_j=W_{2j+2,j+1}$.

\begin{lemma}\label{lem:kernel-rec}
For all $\ell,a\geq0$,
\begin{equation}\label{eq:kernel-rec}
\begin{aligned}
d_{\ell+1}\omega_{\ell+1,a}={}&a_{\ell+a}\,\omega_{\ell,a}+\sum_{j\geq0}(-1)^j\Big(t\,A_j(\ell-j,a+2j+1)\,\omega_{\ell-j,a+2j+1}\\
&+B_j(\ell-j,a+2j+2)\,\omega_{\ell-j,a+2j+2}\\
&+C_j(\ell-j-1,a+2j+2)\,\omega_{\ell-j-1,a+2j+2}\Big),
\end{aligned}
\end{equation}
where $a_i=1-d_i$ is extended by the polynomial $d_i$ for $i>N$ and terms with a negative index are omitted.
\end{lemma}
\begin{proof}
Since $S^aZ^\ell=\sum_i\binom ai z^{\ell+i}w^{\ell+a-i}$, Vandermonde's identity $\binom{\ell+i}k=\sum_{p+q=k}\binom\ell q\binom ip$ and $\binom ai\binom ip=\binom ap\binom{a-p}{i-p}$ give
\[
d(\theta_z)\,S^aZ^\ell=\sum_{p,q\geq0}\nu_{p+q}\binom\ell q\binom ap\,z^pS^{a-p}Z^\ell .
\]
Hence $(\mathcal L_z-\mathcal L_w)S^aZ^\ell$ equals $(z-w)S^aZ^\ell$ plus $\sum_{p,q}\nu_{p+q}\binom\ell q\binom ap S^{a-p}Z^\ell$ times $(z^{p-1}-w^{p-1})+t(z^p-w^p)-(z^{p+1}-w^{p+1})$. Divide by $z-w$, using $\chi_n=(z^{n+1}-w^{n+1})/(z-w)=\sum_{0\leq j\leq\lfloor n/2\rfloor}(-1)^j\binom{n-j}jS^{n-2j}Z^j$ for $n\geq0$, $\chi_{-1}=0$ and $(z^{-1}-w^{-1})/(z-w)=-Z^{-1}$. Since $\sum_q\binom\ell q\nu_{p+q}=f_p(\ell)$, the quotient is
\[
\begin{aligned}
&-d_\ell S^aZ^{\ell-1}+a_{\ell+a}S^aZ^\ell\\
&+\sum_{j\geq0}(-1)^j\Big(tA_j(\ell,a)S^{a-2j-1}Z^{\ell+j}+B_j(\ell,a)S^{a-2j-2}Z^{\ell+j}+C_j(\ell,a)S^{a-2j-2}Z^{\ell+j+1}\Big).
\end{aligned}
\]
Here $-d_\ell S^aZ^{\ell-1}$ comes from $p=0$; the term $a_{\ell+a}S^aZ^\ell$ combines $S^aZ^\ell$ with the parts $j=0$ of $-\chi_p$, which sum to $-d_{\ell+a}S^aZ^\ell$ by Vandermonde; and the three families from $\chi_{p-2}$, $t\chi_{p-1}$ and $-(\chi_p-S^p)$, whose coefficients $\binom{p-2-j}j$, $\binom{p-1-j}j$ and $\binom{p-1-j}{j+1}$ are those of $B_j$, $A_j$ and $C_j$. Summing over $\ell,a$ with weights $\omega_{\ell,a}$ and taking the coefficient of $S^aZ^\ell$ gives \eqref{eq:kernel-rec}, because the monomials $S^aZ^\ell$ are linearly independent.
\end{proof}

\paragraph{The log-concavity inequality.}
\begin{lemma}\label{lem:lc-ineq}
Under the hypotheses of Theorem~\ref{thm:lc}, for $i,j\geq1$ and $\ell\geq0$,
\[
f_i(\ell)\,f_j(\ell+1)\geq d_{\ell+1}\,f_{i+j}(\ell).
\]
\end{lemma}
\begin{proof}
For a nonnegative finitely supported sequence $u$ that is log-concave without internal zeros, $v_i=u_i+u_{i+1}$ has the same property, since
$v_i^2-v_{i-1}v_{i+1}=(u_i^2-u_{i-1}u_{i+1})+(u_{i+1}^2-u_iu_{i+2})+(u_iu_{i+1}-u_{i-1}u_{i+2})$
and the last term is nonnegative by log-concavity. The sequence $(0,\nu_1,\nu_2,\ldots)$ has this property, and $f_i(\ell)$ is obtained from it by applying $u\mapsto v$ $\ell$ times. For such a sequence $g$, balanced products dominate: $g_ig_j\geq g_0g_{i+j}$ and $g_ig_{j+1}\geq g_1g_{i+j}$ for $i,j\geq1$. Adding these with $g=f(\ell)$ and using $g_j+g_{j+1}=f_j(\ell+1)$ and $g_0+g_1=f_0(\ell+1)=d_{\ell+1}$ gives the claim.
\end{proof}

\begin{lemma}\label{lem:detours}
Under the hypotheses of Theorem~\ref{thm:lc}, for all $j\geq1$ and $\ell,a\geq0$,
\[
\begin{aligned}
A_0(\ell,a)B_{j-1}(\ell+1,a-1)+B_0(\ell,a)A_{j-1}(\ell+1,a-2)&\geq4d_{\ell+1}A_j(\ell,a),\\
B_0(\ell,a)B_{j-1}(\ell+1,a-2)&\geq4d_{\ell+1}B_j(\ell,a),\\
B_0(\ell,a)C_{j-1}(\ell+1,a-2)&\geq4d_{\ell+1}C_j(\ell,a),
\end{aligned}
\]
with $W_{p,q}(\ell,a)=0$ for $a<0$.
\end{lemma}
\begin{proof}
Use $A_0(\ell,a)\geq af_1(\ell)$ and $B_0(\ell,a)\geq\binom a2f_2(\ell)$, expand the second factor, apply Lemma~\ref{lem:lc-ineq} to $f_1(\ell)f_n(\ell+1)$ and $f_2(\ell)f_n(\ell+1)$, and use $a\binom{a-1}n=(n+1)\binom a{n+1}$ and $\binom a2\binom{a-2}n=\binom{n+2}2\binom a{n+2}$. After division by $d_{\ell+1}$ the three left sides are bounded below by $\sum_m\gamma_m\binom amf_m(\ell)$ with
\[
\gamma_m=\binom{m+1}2\binom{m-j-2}{j-1},\qquad\binom m2\binom{m-j-3}{j-1},\qquad\binom m2\binom{m-j-2}j,
\]
while the right sides are $4\sum_m\gamma'_m\binom amf_m(\ell)$ with $\gamma'_m=\binom{m-j-1}j$, $\binom{m-j-2}j$ and $\binom{m-j-1}{j+1}$. Where $\gamma'_m>0$ the quotients are
\[
\frac{jm(m+1)}{2(m-j-1)},\qquad\frac{jm(m-1)}{2(m-j-2)},\qquad\frac{(j+1)m(m-1)}{2(m-j-1)},
\]
on $m\geq2j+1$, $m\geq2j+2$ and $m\geq2j+2$ respectively. They increase with $j$ for fixed $m$; at $j=1$ their minima over integers are $5$, $5$ and $6$. All quotients exceed $4$, and $f_m(\ell)\geq0$.
\end{proof}

\begin{proof}[Proof of Theorem~\ref{thm:lc}]
Let $\Omega=\{(\ell,a):\ell,a\geq0,\ \ell+a\leq N\}$, on which $\omega$ is supported. Read \eqref{eq:kernel-rec} as moves into the node $(\ell+1,a)$, each divided by $d_{\ell+1}>0$: the straight move from $(\ell,a)$ with weight $a_{\ell+a}\geq0$ on $\Omega$, and for each $j$ the three moves from $(\ell-j,a+2j+1)$, $(\ell-j,a+2j+2)$, $(\ell-j-1,a+2j+2)$ with weights $tA_j$, $B_j$, $C_j$. Every move raises $\ell$. Let $\mathcal U$ be the straight moves, $\mathcal P$ the moves with $j$ even and $\mathcal Q$ the moves with $j$ odd, all with nonnegative weights, so that $\omega=\beta+(\mathcal U+\mathcal P-\mathcal Q)\omega$ with $\beta=\omega|_{\ell=0}$.

We claim $\mathcal P^2\geq4\mathcal Q$ entrywise. A move of $\mathcal Q$ with index $j$ from $x=(\ell,a)$ is matched by the two-step paths in $\mathcal P^2$ that start with $A_0$ or $B_0$ from $x$ (index $0$, landing at $(\ell+1,a-1)$ or $(\ell+1,a-2)$, which lie in $\Omega$) and continue with the move of the same type and index $j-1$ (even), as in Lemma~\ref{lem:detours}: $A_j$ by $A_0$ then $B_{j-1}$ and $B_0$ then $A_{j-1}$; $B_j$ by $B_0$ then $B_{j-1}$; $C_j$ by $B_0$ then $C_{j-1}$. These paths end at the target of the matched move, carry the same power of $t$, and have the extra denominator $d_{\ell+1}$ of the intermediate node. Distinct moves of $\mathcal Q$ have distinct displacements, so their paths are distinct, and Lemma~\ref{lem:detours} gives the claim.

Now show by induction on the level that $\omega\geq0$ and $\omega\geq\frac12\mathcal P\omega$. Both hold at $\ell=0$. If they hold below a node $y$ at level $\ell+1\geq1$, then, all weights being nonnegative, $(\mathcal P\omega)(y)\geq\frac12(\mathcal P^2\omega)(y)\geq2(\mathcal Q\omega)(y)$, and
\[
\omega(y)=(\mathcal U\omega)(y)+(\mathcal P\omega)(y)-(\mathcal Q\omega)(y)\geq\tfrac12(\mathcal P\omega)(y)\geq0.
\]
Hence all $\omega_{\ell,a}\geq0$, so $K_k$ has nonnegative coefficients in $z,w$, so $L_{mn}^k\geq0$, $H_{mn}^k\geq0$ and $g(m,n;k)\geq0$.
\end{proof}

\begin{remark}\label{rem:bernstein}
The argument is valid coefficientwise in $t$. With independent parameters in $yp_i=a_ip_{i+1}+ud_ip_i+vd_ip_{i-1}$, weighted homogeneity ($\deg y=\deg u=1$, $\deg v=2$) then shows that each linearization coefficient is a polynomial in $(u,v)$ with nonnegative coefficients: for $D=m+n-k\geq0$ it has the form $\sum_{0\leq j\leq\lfloor D/2\rfloor}h_ju^{D-2j}v^j$ with $h_j\geq0$, and it vanishes for $D<0$. Consequently every $g(m,n;k)$ of the array of Theorem~\ref{thm:lc}, as a function of $r$, equals $\sum_jh_jr^j(1-r)^{D-2j}$; multiplying the $j$-th term by $(r+(1-r))^j=1$ shows that its coefficients in the Bernstein basis of degree $D$ are $\binom Ds^{-1}\sum_{j\leq s\leq2j}h_j\binom j{s-j}\geq0$.
\end{remark}

Szwarc~\cite[\S5, Problem~2]{szwarc} records the positive range of the Krawtchouk polynomials as $0<p_0\leq1/2$.\footnote{In the author version of~\cite{szwarc}, \S5, Problem~2, printed p.~125, the displayed weights $\binom Nn p^{N-n}(1-p)^n$ have mean $N(1-p)$, whereas the displayed recurrence gives $K_1(x)=(x-pN)/((1-p)N)$, so that $\int K_1\,d\mu=(1-2p)/(1-p)$. The corollary fixes the convention directly: it uses the displayed weights and the normalization at the point carrying the weight $p^N$.}

\begin{corollary}[Full positive Krawtchouk range]\label{cor:kraw}
For every $N\geq1$ and $0<p_0\leq1/2$, Krawtchouk polynomials orthogonal for binomial weights $\binom Nj p_0^{N-j}(1-p_0)^j$ have nonnegative linearization when normalized at $j=0$.
\end{corollary}
\begin{proof}
To fix the convention without relying on a named recurrence, put $r=p_0/(1-p_0)\in(0,1]$, and let $Y_1,\ldots,Y_N$ be independent variables with values $1,-r$ of probabilities $p_0,1-p_0$. Their mean is zero and $Y_i^2=(1-r)Y_i+r$. The normalized elementary symmetric sums
\[
K_i(j)=\binom Ni^{-1}\sum_{|S|=i}\prod_{l\in S}Y_l
\]
depend only on the number $j$ of negative coordinates. They are orthogonal, because the expectation of two distinct squarefree monomials is zero: some independent mean-zero variable occurs once. Their weights in $j$ are exactly the stated binomial weights. With $x=N^{-1}\sum_lY_l=1-j/(N(1-p_0))$, multiplication by $x$ gives
\[
xK_i=\frac{N-i}{N}K_{i+1}+(1-r)\frac iN K_i+r\frac iN K_{i-1}.
\]
Indeed, in multiplying the elementary sum by $\sum_lY_l$, indices outside a subset give $(i+1)$ copies of each subset of size $i+1$; the linear part of $Y_l^2$ gives $(1-r)i$ copies of each subset of size $i$; and the constant part gives $r(N-i+1)$ copies of each subset of size $i-1$. Division by the binomial normalizers gives the displayed recurrence, including $i=N$. Its leading coefficient in $x$ is positive, its degree is $i$, and its value at $j=0$ is one. These properties identify the normalized orthogonal polynomials for the binomial measure. This is the array $K_{N,r}$; its profile $d_i=i/N$ has Newton coefficients $(1/N,0,\ldots,0)$, so Theorem~\ref{thm:lc} applies.
The range is sharp in this endpoint normalization: if $1/2<p_0<1$, then
$r=p_0/(1-p_0)>1$ and $g(1,1;1)=b_1=(1-r)/N<0$.
For $N=1$, this coefficient follows from the terminal relation.
\end{proof}

\begin{corollary}[Geometric arrays]\label{cor:additional}
For every $N\geq2$, $q>1$, and $0<r\leq1$, the array $G_{N,q,r}$ of \eqref{eq:parameters} has nonnegative linearization. In the normalization $R_0=1$, $R_1=x$, it is outside the Krawtchouk recurrence class.
\end{corollary}
\begin{proof}
Its increments are $d_i-d_{i-1}=q^{i-1}/[N]_q$, so its Newton coefficients $\nu_k=(q-1)^{k-1}/[N]_q$ are log-concave without internal zeros and Theorem~\ref{thm:lc} applies. For a Krawtchouk recurrence in this normalization the increments of $b_i+c_i$ are constant. Here their ratio is $q>1$ and $N\geq2$, so they are not constant. An affine change of the support transporting the same endpoint-normalized basis leaves the multiplication coefficients in the basis unchanged; it cannot remove this distinction.
\end{proof}

The arrays $G_{N,q,r}$ are the normalized recurrences of Koelink's construction~\cite[Theorems~3.1(ii) and~7.1]{koelink} with $r=1/p$; their polynomial interpretation is the dual $q$-Krawtchouk interpretation described after his Theorem~7.8. A proof of Corollary~\ref{cor:additional} by a finite algebra with positive multiplication, which needs only the geometric increments, is given in~\cite{P3}.

For some parameters both families come from distance-regular graphs. Dividing the intersection array of such a graph by its valency gives an array satisfying \eqref{eq:basic}, and its linearization coefficients are positive multiples of the intersection numbers $p^k_{mn}$, hence nonnegative~\cite{bcn,vdkt}. For an integer $s\geq2$, the Hamming graph $H(N,s)$ gives $K_{N,r}$ with $r=1/(s-1)$; a dual polar graph of diameter $N$ with classical parameters $(N,q,0,q^e)$ gives $G_{N,q,r}$ with $r=q^{-e}$. Corollaries~\ref{cor:kraw} and~\ref{cor:additional} cover the full parameter ranges.

\section{Szwarc's quadratic system}\label{sec:szwarc}

Problem~2 of~\cite[\S5]{szwarc} does not stop at the Krawtchouk recurrence. With his parameter written as $p_0$, Szwarc introduces the symmetric recurrence of length $2N+1$
\begin{equation}\label{eq:szwarc-q}
yQ_{2n}=(N-n)Q_{2n+1}+nQ_{2n-1},\qquad
yQ_{2n-1}=(1-p_0)Q_{2n}+p_0Q_{2n-2},
\end{equation}
records that $Q_{2n}(y)$ is the Krawtchouk polynomial $K_n(y^2)$ of his recurrence, states that the $Q_n$ have nonnegative linearization for $p_0\leq1/2$, and asks for a general criterion that would handle this case. His Theorem~3.3 for infinite symmetric systems does not apply. In the normalization \eqref{eq:szwarc-q} its monotonicity hypotheses hold, and its hypothesis $\alpha_n\leq\gamma_n$ holds at odd indices exactly when $p_0\leq1/2$; at an even index $2n$ it reads $n\leq N-n$ and fails in the upper half. Some failure of this kind is forced in every finite system, whose terminal upward coefficient vanishes. We derive the positivity of \eqref{eq:szwarc-q} from Theorem~\ref{thm:lc} through a quadratic decomposition of symmetric finite systems.

Put $M=2N$, $R_{2n}(t)=Q_{2n}(\sqrt Nt)$ and $R_{2n-1}(t)=\sqrt N\,Q_{2n-1}(\sqrt Nt)$. Substitution in \eqref{eq:szwarc-q} shows that the $R_k$ satisfy \eqref{eq:recurrence}, with $N$ replaced by $M$, for the array
\begin{equation}\label{eq:qarray}
b_k=0,\qquad
(a_{2n},c_{2n})=\Big(\frac{N-n}{N},\frac nN\Big),\qquad
(a_{2n-1},c_{2n-1})=(1-p_0,\,p_0).
\end{equation}
It satisfies \eqref{eq:basic} with $M$ in place of $N$ whenever $0<p_0<1$. The two normalizations differ by positive factors and a positive dilation, so they have the same product coefficients up to positive factors.

Szwarc uses the quadratic transformation~\cite[\S3.5]{szwarc} to pass from a symmetric system to its even part. In~\cite{szwarc1992} a parity-split monotonicity criterion for symmetric infinite systems gives property (P) of the generalized Chebyshev polynomials and, through their even parts, Gasper's result for the Jacobi polynomials with $\alpha\geq\beta$ and $\alpha+\beta+1\geq0$; in~\cite{szwarc1996} the same route gives property (P) of Askey--Wilson polynomials for a range of parameters. In both, positivity passes from a symmetric system to its even part. For the finite systems of classical association schemes, bipartite halves and antipodal quotients give quadratic transformations of the polynomials, which Chihara and Stanton~\cite{chiharastanton} derive from a quadratic transformation of the Askey--Wilson polynomials. The decomposition of a symmetric system into its even and odd components, and the identities relating their products, are standard, and we claim no priority for them. What we use is the converse implication for finite systems: positivity of a symmetric system, including the products beyond the terminal degree, follows from positivity of its two components, both read off from its coefficients, which need not come from a scheme. Any criterion for nonsymmetric finite systems therefore becomes a criterion for symmetric ones. For a family of infinite systems a result of this kind is known: Kahler~\cite{kahler} solved Problem~1 of~\cite[\S5]{szwarc} by showing that the generalized Chebyshev polynomials, the quadratic transforms of the Jacobi polynomials, have property (P) exactly for the parameters for which the Jacobi polynomials $R_n^{(\alpha,\beta)}$ have it. His result concerns that family; Theorem~\ref{thm:decomp} below is a reduction for arbitrary finite symmetric systems.

\begin{theorem}[Quadratic decomposition]\label{thm:decomp}
Let $N\geq2$, and let $a_k,c_k$ $(0\leq k\leq2N)$ with $b\equiv0$ satisfy \eqref{eq:basic} with $2N$ in place of $N$. Put $a_k=c_k=0$ for $k\notin[0,2N]$, and define the even component
\[
\hat a_n=a_{2n}a_{2n+1},\quad \hat b_n=c_{2n}a_{2n-1}+a_{2n}c_{2n+1},\quad \hat c_n=c_{2n}c_{2n-1}\qquad(0\leq n\leq N)
\]
and the odd component
\[
\check a_n=a_{2n+1}a_{2n+2},\quad \check b_n=a_{2n+1}c_{2n+2}+c_{2n+1}a_{2n},\quad \check c_n=c_{2n+1}c_{2n}\qquad(0\leq n<N).
\]
After the affine normalization $(\hat a_n,\hat b_n,\hat c_n)\mapsto(\hat a_n,\hat b_n-\hat b_0,\hat c_n)/\hat a_0$, and likewise for the odd component, these arrays satisfy \eqref{eq:basic} on the indices $0,\ldots,N$ and $0,\ldots,N-1$. If both have nonnegative linearization, so does the symmetric array. More precisely, with $\hat g$ and $\check g$ the product coefficients of the components and terms with indices out of range omitted,
\[
g(2a,2b;2c)=\hat g(a,b;c),\qquad
g(2a,2b+1;2c+1)=a_{2a}\,\check g(a,b;c)+c_{2a}\,\check g(a-1,b;c),
\]
the coefficients with $m+n+k$ odd vanish, and $g(2a+1,2b+1;2c)\int R_{2c}^2\,d\mu=g(2c,2a+1;2b+1)\int R_{2b+1}^2\,d\mu$. Conversely, nonnegative linearization of the symmetric array implies that of its even component.
\end{theorem}
\begin{proof}
Since $b\equiv0$, the matrix $J$ is conjugate to $-J$ by $\operatorname{diag}((-1)^i)$, which fixes $e_0$; hence the measure of Lemma~\ref{lem:finite} is symmetric and $R_k(-t)=(-1)^kR_k(t)$. By symmetry $R_{2n}(t)=P_n(t^2)$ and $R_{2n+1}(t)=tL_n(t^2)$ for polynomials $P_n$, $L_n$ of degree $n$. Applying \eqref{eq:recurrence} twice gives, with $x=t^2$,
\[
xP_n=\hat a_nP_{n+1}+\hat b_nP_n+\hat c_nP_{n-1},\qquad
xL_n=\check a_nL_{n+1}+\check b_nL_n+\check c_nL_{n-1}.
\]
The row sums of the components equal $a_{2n}+c_{2n}=1$ and $a_{2n+1}+c_{2n+1}=1$; the coefficients $\hat a_n,\check a_n$ are positive below the last index and vanish there, and $\hat c_n,\check c_n>0$ for $n\geq1$. Hence the normalized arrays, which describe the same polynomials in the variable $(x-\hat b_0)/\hat a_0$ (respectively $(x-\check b_0)/\check a_0$), satisfy \eqref{eq:basic}. The terminal relation $tR_{2N}=c_{2N}R_{2N-1}$ on the support $\Sigma$ of the symmetric system, applied once and twice, gives the terminal relations of $P$ on $\{t^2:t\in\Sigma\}$ and of $L$ on $\{t^2:t\in\Sigma,\ t\neq0\}$. These sets have $N+1$ and $N$ points, because $\Sigma$ is symmetric with $2N+1$ points and contains $0$. By Lemma~\ref{lem:finite} they are the supports of the two components.

The identity $tR_{2n}=a_{2n}R_{2n+1}+c_{2n}R_{2n-1}$ gives $P_n=a_{2n}L_n+c_{2n}L_{n-1}$ on the support of $L$; for $n=N$ it is the terminal relation, with $a_{2N}=0$. Hence, on $\Sigma\setminus\{0\}$,
\[
R_{2a}R_{2b+1}=t\big(a_{2a}L_aL_b+c_{2a}L_{a-1}L_b\big)=\sum_c\big(a_{2a}\check g(a,b;c)+c_{2a}\check g(a-1,b;c)\big)R_{2c+1},
\]
and both sides vanish at $t=0$. Similarly $R_{2a}R_{2b}=(P_aP_b)(t^2)=\sum_c\hat g(a,b;c)R_{2c}$ on $\Sigma$. Uniqueness of expansions on $\Sigma$ gives the displayed formulas. The vanishing for $m+n+k$ odd and the relation for two odd indices follow from the symmetry of the measure and from $g(m,n;k)\int R_k^2\,d\mu=\int R_mR_nR_k\,d\mu$. The converse holds because $\hat g(a,b;c)=g(2a,2b;2c)$.
\end{proof}

Combined with Theorem~\ref{thm:lc}, the decomposition gives a parity-split criterion of the same shape as Szwarc's Theorem~3.3: separate conditions on the even-indexed and the odd-indexed coefficients.

\begin{corollary}[Parity-split criterion]\label{cor:parity}
Let $N\geq1$ and let a symmetric array $(b\equiv0)$ on the indices $0,\ldots,2N$ satisfy \eqref{eq:basic}. Suppose that
\begin{enumerate}
\item the odd downward coefficients are constant, $c_{2n+1}=p_0$, with $0<p_0\leq1/2$;
\item the even downward coefficients $d_n=c_{2n}$ have positive increments $d_n-d_{n-1}$, $1\leq n\leq N$, forming a geometric progression with ratio $q\geq1$.
\end{enumerate}
Then the array has nonnegative linearization. With $r=p_0/(1-p_0)$ and $N\geq2$, its normalized even and odd components are \eqref{eq:parameters} with parameters $(N,q,r)$ and $(N-1,q,r/q)$, respectively. For $N=1$ the odd component consists only of $L_0=1$ and the array is \eqref{eq:qarray}. If $1/2<p_0<1$ in (i) and (ii) holds, then $g(2,2;2)<0$.
\end{corollary}
\begin{proof}
For $N=1$ the array is $a_0=1$, $c_1=p_0$, $a_1=1-p_0$, $c_2=1$, with support $\{-1,0,1\}$ and $R_1=t$, $R_2=(t^2-p_0)/(1-p_0)$. On the support $R_1^2=(1-p_0)R_2+p_0$, $R_1R_2=R_1$ and $R_2^2=\frac{p_0}{1-p_0}+\frac{1-2p_0}{1-p_0}R_2$, as evaluation at $t=0,\pm1$ shows; all coefficients are nonnegative for $p_0\leq1/2$, and $g(2,2;2)<0$ for $p_0>1/2$. Let $N\geq2$. The even component has $\hat a_n=(1-d_n)(1-p_0)$, $\hat b_n=d_n(1-p_0)+(1-d_n)p_0$ and $\hat c_n=d_np_0$, so $\hat a_0=1-p_0$, $\hat b_0=p_0$, and its normalization is $(1-d_n,\,(1-r)d_n,\,rd_n)$. Its profile $d_n$ has geometric increments, so its Newton coefficients are proportional to $(q-1)^{k-1}$ and log-concave without internal zeros, and $0<r\leq1$. For the odd component, $\check a_n=(1-p_0)(1-d_{n+1})$, $\check b_n=(1-p_0)d_{n+1}+p_0(1-d_n)$ and $\check c_n=p_0d_n$. Geometric increments give $d_{n+1}-d_1=qd_n$, so the normalization is
\[
\Big(\frac{1-d_{n+1}}{1-d_1},\ \frac{(q-r)d_n}{1-d_1},\ \frac{rd_n}{1-d_1}\Big).
\]
Its profile $qd_n/(1-d_1)$ has geometric increments with ratio $q$, and its downward coefficients are $r/q\leq1$ times the profile. Theorem~\ref{thm:lc} applies to both components, and Theorem~\ref{thm:decomp} concludes. Comparison with \eqref{eq:parameters}, using $d_n=[n]_q/[N]_q$ and $1-d_1=q[N-1]_q/[N]_q$, identifies the parameters. For $p_0>1/2$, $g(2,2;2)=\hat g(1,1;1)$ is the normalized diagonal coefficient $(1-r)d_1<0$ (for $N=1$ see above).
\end{proof}

\begin{theorem}[Szwarc's quadratic system]\label{thm:szwarc}
Let $N\geq1$ and $0<p_0\leq1/2$. The array \eqref{eq:qarray} has a symmetric positive orthogonality measure on $2N+1$ points, and for all $0\leq m,n\leq 2N$
\[
R_mR_n=\sum_{k=0}^{2N}g(m,n;k)R_k,\qquad g(m,n;k)\geq0,\qquad\sum_k g(m,n;k)=1
\]
on that support. If $1/2<p_0<1$, then $g(2,2;2)<0$.
\end{theorem}
\begin{proof}
Lemma~\ref{lem:finite} gives a positive measure on $2N+1$ points, and it is symmetric because $b\equiv0$ (proof of Theorem~\ref{thm:decomp}). The array \eqref{eq:qarray} satisfies the hypotheses of Corollary~\ref{cor:parity} with $c_{2n+1}=p_0$ and $d_n=n/N$, that is, with $q=1$, which gives nonnegative linearization for $p_0\leq1/2$ and $g(2,2;2)<0$ for $1/2<p_0<1$. The coefficients sum to one because $R_k(1)=1$ and $1$ lies in the support, since $J\boldsymbol1=\boldsymbol1$.
\end{proof}

The even component of \eqref{eq:qarray} is the Krawtchouk array $K_{N,r}$ with $r=p_0/(1-p_0)$, so that $Q_{2n}(y)$ is a Krawtchouk polynomial in $y^2$, as Szwarc notes, and its odd component is $K_{N-1,r}$. A second proof of Theorem~\ref{thm:szwarc}, by a graded version of the algebra of~\cite{P3}, and explicit coefficient neighborhoods of \eqref{eq:qarray} are given there.

\section{Limits of the criterion}\label{sec:limits}

\paragraph{No uniform margin.} A criterion valid for all $N$ cannot leave a margin of fixed size around the Krawtchouk arrays. Let $\rho_N(r)$ be the supremum of the radii $\rho$ such that every array satisfying \eqref{eq:basic}, with $b_i\geq0$ and with coefficients that differ from those of $K_{N,r}$ by at most $\rho$, has property (P). Raising the top mass of the orthogonality measure of $K_{N,r}$ (Proposition~\ref{prop:raise}) gives $\rho_N(r)\leq3(1-r^2)^N$ whenever $(1-r^2)^N\leq\frac12$ and $3(1-r^2)^N\leq(1-r)/N$ (Corollary~\ref{cor:radius-exp}); both are proved in Appendix~\ref{app:radius}. So, for every fixed $r\in(0,1)$, no set of arrays with property (P) contains coefficient balls of a radius independent of $N$ around all $K_{N,r}$, and the log-concavity condition places these arrays on its edge. This is an upper bound for fixed $r$ and $N\to\infty$ only. Lower bounds, the exact order $(1-r)^N$ for fixed $N$ as $r\to1$, and the constraints at $r=1$ are proved in~\cite{P3}.

\paragraph{How sharp is log-concavity?}
The Krawtchouk arrays, with only $\nu_1\neq0$, and the arrays $G_{N,q,r}$, with equality in every log-concavity inequality, both lie on the boundary of the hypothesis of Theorem~\ref{thm:lc}. Positivity of the Newton coefficients of the profile does not suffice. For $q>1$ and $0<\lambda<1$ the combination $(1-\lambda)K_{N,r}+\lambda G_{N,q,r}$ is a stochastic array with constant ratio $b_i/c_i$ and positive Newton coefficients, and $\nu_2^2<\nu_1\nu_3$ for $N\geq3$ (Section~\ref{sec:lc}). Exact computation, machine-verified in Lean, gives $g(N,N;N)<0$ for seven such combinations, with $r$ between $2/7$ and $4/5$ (Appendix~\ref{app:search}). At these parameters no sufficient criterion whose admitted set is convex, such as a system of linear inequalities in the coefficients, can admit both $K_{N,r}$ and $G_{N,q,r}$.

Log-concavity is not necessary either. For $0<\mu<1$ and $q>1$ let $D(i)=(1-\mu)i+\mu[i]_q$, whose Newton coefficients are $\nu_1=1$ and $\nu_k=\mu(q-1)^{k-1}$ for $k\geq2$; thus $\nu_1\nu_3/\nu_2^2=1/\mu$, and the other log-concavity inequalities hold with equality. For $\mu=\frac1{10}$ and $q=\frac{11}{10}$, the constant-ratio arrays with profile $d_i=D(i)/D(N)$ have property (P) for every $r\in(0,1]$ and $N\in\{6,8,10,12\}$ (machine-verified), whereas $\mu=\frac1{100}$, the same $q$ and $N=6$ give $g(6,6;6)<0$ at $r=1$. Let
\[
V_N=\inf\{\mu^{-1}:0<\mu<1,\ q>1,\ \text{the array with profile }D/D(N)\text{ and }r=1\text{ fails (P)}\},
\]
with $\inf\varnothing=\infty$. Numerical searches for the smallest failing $1/\mu$ are recorded in Appendix~\ref{app:search}.

\section{Open problems}\label{sec:open}

\paragraph{The exact boundary.} Determine the boundary between log-concavity and positivity of the Newton coefficients of the profile: the set of profiles for which every constant-ratio array has property (P). In particular, decide whether the threshold $V_N$ of Section~\ref{sec:limits} tends to $1$ as $N\to\infty$. The examples of Section~\ref{sec:limits} suggest that when log-concavity fails, the failures of (P) are driven by an acceleration of the increments near the terminal index.

\paragraph{The stability radius.} For fixed $0<r<1$, determine the order of $\rho_N(r)$ as $N\to\infty$, and in particular decide whether $\liminf_N\rho_N(r)^{1/N}$ is positive. Corollary~\ref{cor:radius-exp} gives $\limsup_N\rho_N(r)^{1/N}\leq1-r^2$.

\paragraph{Other ratios.} Theorem~\ref{thm:lc} requires $b_i/c_i$ to be constant. The companion paper~\cite{P2} admits ratios $c_i/d_i$ that increase with the index. On log-concave Newton profiles with $b$ nondecreasing, it admits linear and convex quadratic ratios, ratios of any degree whose Newton coefficients are log-concave without internal zeros, and cubic ratios under a single inequality. Criteria for other nonconstant ratios are open.

\paragraph{Monotonicity-type criteria.} Szwarc's Theorems~3.2 and~3.3 in~\cite{szwarc} are inequality conditions of monotonicity type for infinite systems. By Theorem~\ref{thm:decomp}, any such criterion for nonsymmetric finite systems, valid including the terminal indices, gives one for symmetric finite systems. Inequality conditions of this type for all products except those whose pairwise index sums all exceed $N$ are given in~\cite{P3}; for the remaining, deep products, Theorem~\ref{thm:lc} and its extensions in~\cite{P2} give shape conditions that control them.

\appendix
\section{The upper bound for the stability radius}\label{app:radius}

This appendix proves Corollary~\ref{cor:radius-exp}. Let
\[
\operatorname{dist}(\widetilde a,\widetilde b,\widetilde c)=\max_{0\leq i\leq N}\max\Bigl(\Bigl|\widetilde a_i-\frac{N-i}N\Bigr|,\Bigl|\widetilde b_i-\frac{(1-r)i}N\Bigr|,\Bigl|\widetilde c_i-\frac{ri}N\Bigr|\Bigr)
\]
be the coefficient distance to $K_{N,r}$.

\begin{lemma}[Corner value]\label{lem:corner-value}
For $0<r\leq1$, the array $K_{N,r}$ has $g(N,N;N)=(1-r)^N$.
\end{lemma}
\begin{proof}
Use the model in the proof of Corollary~\ref{cor:kraw}, with $p_0=r/(1+r)$. There $R_N=K_N=Y_1\cdots Y_N$ and $K_k=\binom Nk^{-1}\sum_{|S|=k}Y_S$, where $Y_S=\prod_{l\in S}Y_l$. Since $Y_l^2=(1-r)Y_l+r$,
\[
R_N^2=\prod_{l=1}^N\bigl((1-r)Y_l+r\bigr)=\sum_S(1-r)^{|S|}r^{N-|S|}Y_S.
\]
Both sides of $R_N^2=\sum_kg(N,N;k)R_k$ depend only on the number $j$ of negative coordinates, and the identity holds on the support, that is, for every $j$; so it holds as an identity of functions of $(Y_1,\ldots,Y_N)$. The $2^N$ monomials $Y_S$ are orthogonal for the product measure, with $\mathbb E[Y_S^2]=r^{|S|}>0$, hence linearly independent. The monomial $Y_{\{1,\ldots,N\}}$ occurs only in $K_N$, with coefficient one, which gives $g(N,N;N)=(1-r)^N$.
\end{proof}

Thus the corner coefficient vanishes at $r=1$. The next proposition raises the top mass of the orthogonality measure.

\begin{proposition}[Raising the top mass]\label{prop:raise}
Let $N\geq1$ and $0<r<1$. The measure of $K_{N,r}$ has nodes $x_n=1-(1+r)n/N$ and masses $w_n=\binom Nnr^{N-n}/(1+r)^N$. For $\tau>0$ let $\mu^\tau$ be this measure with the mass at $x_0=1$ raised to $(1+\tau)w_0$. Then $\mu^\tau$ has
\[
g(N,N;N)=\frac{(1+\tau)^2(I-w_0)+w_0}{(1+\tau)(h-w_0)+w_0},\qquad I=\bigl(r(1-r)\bigr)^N,\quad h=r^N,
\]
which is negative exactly when $\tau>\tau_N^*=\bigl(1-(1-r^2)^N\bigr)^{-1/2}-1$. After the affine normalization $R_n(1)=1$, $R_1=x$, its recurrence coefficients differ from those of $K_{N,r}$ by at most $3\tau\max_nw_n$.
\end{proposition}
\begin{proof}
Write $\pi_n$ and $S_n=\sum_{j<n}\pi_j$ for $K_{N,r}$, and put $\varepsilon=\tau w_0$. The polynomial $(1+\varepsilon S_n)R_n-\varepsilon\sum_{j<n}\pi_jR_j$ has degree $n$ and equals $1$ at $x_0$. For $\deg p<n$, its integral against $p$ with respect to $\mu^\tau=\mu+\varepsilon\delta_{x_0}$ is $-\varepsilon p(x_0)+\varepsilon p(x_0)=0$, because $\sum_{j<n}\pi_jR_j$ reproduces $p(x_0)$. Hence it is the polynomial $R_n^\tau$ of $\mu^\tau$.

On the support, dual orthogonality gives $\sum_{j\leq N}\pi_jR_j(x)=\boldsymbol1_{x=x_0}/w_0$, and in particular $S_N+\pi_N=1/w_0$. Therefore $R_N^\tau=(1+\tau)R_N-\tau\boldsymbol1_{x_0}$ on the support. By Lemma~\ref{lem:corner-value} and $\pi_N=r^{-N}$ (Lemma~\ref{lem:finite}), $\int R_N^3\,d\mu=I$; also $\int R_N^2\,d\mu=h$. Hence
\begin{gather*}
\int(R_N^\tau)^3d\mu^\tau=(1+\tau)^3(I-w_0)+(1+\tau)w_0,\\
\int(R_N^\tau)^2d\mu^\tau=(1+\tau)^2(h-w_0)+(1+\tau)w_0,
\end{gather*}
and their ratio is the stated formula. Since $r>r/(1+r)>r(1-r)$, we have $h>w_0>I$. So the denominator is positive, and the numerator is negative exactly when $(1+\tau)^2>w_0/(w_0-I)=\bigl(1-(1-r^2)^N\bigr)^{-1}$.

For the coefficients, comparing leading coefficients and using detailed balance, $\pi_na_n=\pi_{n+1}c_{n+1}$, for both measures gives
\[
a_n^\tau=a_n\frac{1+\varepsilon S_n}{1+\varepsilon S_{n+1}},\qquad c_n^\tau=c_n\frac{1+\varepsilon S_{n+1}}{1+\varepsilon S_n},
\]
and the row sums stay $1$. Since $\varepsilon\pi_n=\tau w_n$, each of $a_n$, $b_n$, $c_n$ changes by at most $\tau w_n$. Put $A_n=1+\varepsilon S_n$ and $t_n=\varepsilon\pi_n=\tau w_n$. Since $b_0^\tau=\varepsilon/(1+\varepsilon)$, the affine coordinate is $y=(1+\varepsilon)x-\varepsilon$, and the normalized coefficients are $\widetilde a_n=(1+\varepsilon)a_n^\tau$, $\widetilde b_n=(1+\varepsilon)b_n^\tau-\varepsilon$, $\widetilde c_n=(1+\varepsilon)c_n^\tau$. For $n\geq1$ we have $A_n\geq1+\varepsilon$, and
\begin{align*}
\widetilde a_n-a_n&=\varepsilon a_n-(1+\varepsilon)a_nt_n/A_{n+1},\\
\widetilde c_n-c_n&=\varepsilon c_n+(1+\varepsilon)c_nt_n/A_n,\\
\widetilde b_n-b_n&=-\varepsilon(a_n+c_n)+(1+\varepsilon)t_n\bigl(a_n/A_{n+1}-c_n/A_n\bigr).
\end{align*}
Because $a_n,c_n\geq0$ and $a_n+c_n\leq1$, each absolute change is at most $\varepsilon+t_n\leq2\tau\max_jw_j$. Row zero is restored exactly. This proves, in particular, the stated bound $3\tau\max_jw_j$ for every $\tau>0$.
\end{proof}

\begin{corollary}[Exponential upper bound for the stability radius]\label{cor:radius-exp}
Let $0<r<1$, and let $\rho_N(r)$ be the supremum of the $\rho$ such that every array satisfying \eqref{eq:basic}, with $\widetilde b_i\geq0$ for all $i$ and with $\operatorname{dist}(\widetilde a,\widetilde b,\widetilde c)\leq\rho$, has property (P). If $(1-r^2)^N\leq\tfrac12$ and $3(1-r^2)^N\leq(1-r)/N$, then
\[
\rho_N(r)\leq3\tau_N^*\max_nw_n\leq3(1-r^2)^N.
\]
In particular $\rho_N(r)=O\bigl((1-r^2)^N\bigr)$ as $N\to\infty$, for every fixed $r\in(0,1)$.
\end{corollary}
\begin{proof}
Put $x=(1-r^2)^N$ and $W=\max_nw_n$. For $0<x\leq\tfrac12$,
\[
(1+x)^2(1-x)-1=x(1-x-x^2)>0,
\]
so $\tau_N^*<x$. Choose $\tau\in(\tau_N^*,x)$. Proposition~\ref{prop:raise} gives a negative corner coefficient and full-array distance at most $3\tau W<3x\leq(1-r)/N$. Thus $\widetilde b_n>0$ for $n\geq1$, while $\widetilde b_0=0$. The exact coefficient formulas and the positive affine scale give $\widetilde a_n>0$ for $n<N$, $\widetilde c_n>0$ for $n>0$, and all the basic boundary and row-sum equalities. Letting $\tau\downarrow\tau_N^*$ proves $\rho_N(r)\leq3\tau_N^*W\leq3x$. For each fixed $r\in(0,1)$, both smallness conditions hold for all sufficiently large $N$, since $N(1-r^2)^N\to0$.
\end{proof}

Thus, for every fixed $r\in(0,1)$, no set of arrays with property (P) contains coefficient balls of a positive radius independent of $N$ around all $K_{N,r}$. Corollary~\ref{cor:radius-exp} bounds $\rho_N(r)$ from above only; it does not determine its order. For each fixed $N$, $\rho_N(r)>0$ for $0<r<1$: explicit perturbation inequalities give $\rho_N(r)\geq\frac12(\min(r,1-r)/(32N))^N$, a lower bound that decays faster than every exponential in $N$, and, for fixed $N\geq2$, $\rho_N(r)$ is of exact order $(1-r)^N$ as $r\to1$; at the symmetric endpoint $r=1$, criteria that admit polynomial deformation rules, or analytic rules of the kind specified there, must satisfy two first-order equalities. These results are proved in~\cite{P3}.

\section{Examples and searches for Section~\ref{sec:limits}}\label{app:search}

\paragraph{Convex combinations.} The combination $(1-\lambda)K_{N,r}+\lambda G_{N,q,r}$ has increments proportional to $1+t\,q^{i-N}$, where $t=\lambda N(q-1)q^{N-1}/((1-\lambda)(q^N-1))$. Exact computation, machine-verified in Lean, gives $g(N,N;N)<0$ at the parameters
\[\begin{gathered}(r,N,q,t)=(4/5,10,10,1),\ (3/5,20,10,2),\ (1/2,20,10,1),\ (2/5,60,10,2),\\ (1/3,100,10,2),\ (3/10,200,2,5),\ (2/7,200,10,5),\end{gathered}\]
while $K_{N,r}$ and $G_{N,q,r}$ have property (P) by Corollaries~\ref{cor:kraw} and~\ref{cor:additional}.

\paragraph{The profiles $D/D(N)$.} For $\mu=\frac1{10}$ and $q=\frac{11}{10}$, the exact computation shows that every coefficient $g(m,n;k)$, $N\in\{6,8,10,12\}$, is a polynomial in $r$ with nonnegative Bernstein coefficients, which gives property (P) for every $r\in(0,1]$.

\paragraph{Searches.} At $r=1$ we searched coarse grids of $\mu$ and $q$, with float64 screening and high-precision confirmation for $N\leq20$ and high-precision arithmetic for $20\leq N\leq80$. The smallest factor $1/\mu$ for which (P) failed was $10^{5/2}$ at $N=4$, $10^2$ for $6\leq N\leq30$, and $10^{3/2}$ at $N=40,60,80$; scripts and outputs are in the companion code~\cite{code}. These are the smallest failures found, not proven thresholds: the observed failures lay in narrow windows of $(\mu,q)$ on the tested grids, and none of these searches is formalized.

\section*{AI-use disclosure}

Claude Opus 5.5 (Anthropic), run in the Claude Code agent; GPT-6 Astra or GPT-6.1 Sol (OpenAI; which of the two was not recorded), run in the Codex agent; and GPT-6 Pro (OpenAI), used in ChatGPT, were used to develop the proofs, the exact checks, the Lean formalization and the text. GPT-6 Pro also reviewed the results. This AI review is not human peer review, and no independent human verification is claimed. The author set the research objective and is responsible for the mathematical claims. The supplement~\cite{supp} gives the full disclosure, including which steps arose with which model.

\section*{Code and formalization}

All theorems, corollaries and propositions, and Lemmas~\ref{lem:finite} and~\ref{lem:corner-value}, are formalized in Lean~4 with Mathlib in the companion code~\cite{code}, whose axiom audit admits only \path{propext}, \path{Classical.choice} and \path{Quot.sound}; ``machine-verified'' refers only to these Lean theorems. The supplement~\cite{supp} lists the Lean declaration for each result and remark, describes the Lean routes that differ from the written proofs, and lists what is not formalized: the comparisons with the literature and the numerical searches concerning $V_N$. The documentation~\cite{docs} describes the exact-arithmetic diagnostics.

All theorems, corollaries and propositions, as well as Lemmas~\ref{lem:finite} and~\ref{lem:corner-value}, are formalized in Lean~4 with Mathlib; the lemmas of Section~\ref{sec:lc} are formalized within the proof of Theorem~\ref{thm:lc}, whose Lean proof follows the written one. Some Lean proofs take another route, through the Lean development of~\cite{P3}. Lean derives the positivity statements of Corollaries~\ref{cor:kraw}, \ref{cor:additional} and~\ref{cor:parity}, and through Corollary~\ref{cor:parity} that of Theorem~\ref{thm:szwarc}, from the geometric-increment criterion of~\cite{P3} rather than from Theorem~\ref{thm:lc}. It proves Lemma~\ref{lem:corner-value}, on which Proposition~\ref{prop:raise} and Corollary~\ref{cor:radius-exp} rest, in the algebra of~\cite{P3}. These Lean proofs verify the statements given here by routes different from the written proofs; the supplement~\cite{supp} gives the details.

The three Lean developments share one library. The substantive dependencies between the papers are these:
\begin{center}
\begin{tabular}{@{}lll@{}}
Paper & Written proofs use & Lean proofs use\\\hline
this paper & --- & \cite{P3} for Corollaries~\ref{cor:kraw}, \ref{cor:additional}, \ref{cor:parity},\\
 & & \quad Theorem~\ref{thm:szwarc} and Lemma~\ref{lem:corner-value}\\
\cite{P2} & this paper & this paper\\
\cite{P3} & --- & this paper for two proofs
\end{tabular}
\end{center}

\end{document}